\documentclass[12pt,reqno]{amsart} 
\usepackage{amsmath, amsthm,
  amssymb,color,graphicx,bm,psfrag,float,mathpazo, verbatim, enumitem,
  fullpage, bm}
\numberwithin{equation}{section}
 
\newcommand{\field}{\mathbb{F}}
\newcommand{\chr}{\text{char}} 
\newcommand{\proj}{\mathbb{P}}
\newcommand{\Sym}{\mathfrak{S}}
\newcommand{\affine}{\mathbb{A}}
\newcommand{\QQ}{\mathbb{Q}}
\newcommand{\ZZ}{\mathbb{Z}}
\newcommand{\RR}{\mathbb{R}}
 
\newcommand{\sig}{\text{sig}}
\newcommand{\steiner}{\text{steiner}}
\newcommand{\leisen}{\text{leisen}} 
\newcommand{\Gal}{\text{Gal}} 
\newcommand{\conic}{\mathcal{C}}
\newcommand{\complex}{\mathbb{C}}
\renewcommand{\mod}[1]{\; \; (\text{mod} \, #1)}
\newcommand{\artin}[3]{\left[ \frac{#1/#2}{#3} \right]}
\newtheorem{Theorem}{Theorem}[section]
\newtheorem{Lemma}[Theorem]{Lemma}
\newtheorem{Proposition}[Theorem]{Proposition}
\newtheorem{Conjecture}[Theorem]{Conjecture} 
\newtheorem{Definition}[Theorem]{Definition} 

\begin{document} 

\title{On the dynamics of the Pappus-Steiner map} 
\author{Jaydeep Chipalkatti \\ (with an appendix by Attila D{\'e}nes)}
\maketitle

\parbox{16.5cm}{ \small
{\sc Abstract:} We  extract a two-dimensional dynamical system from
the theorems of Pappus and Steiner in classical projective
geometry. We calculate an explicit formula for this system, and study its
elementary geometric properties. Then we use Artin reciprocity to characterise 
all sufficiently large primes $p$ for which this system admits periodic points of orders $3$ and $4$ over the field 
$\field_p$; this leads to an unexpected Galois-theoretic conjecture
for $n$-periodic points. We also give a short discussion of
Leisenring's theorem, and show that it leads to the same dynamical
system as the Pappus-Steiner theorem. The appendix contains a
computer-aided analysis of this system over the field of real numbers.} 

\bigskip \bigskip


\medskip 

\tableofcontents

\section{Introduction} 
This article is a study of a two-dimensional dynamical system arising
out of the Pappus-Steiner theorem in classical projective
geometry. The construction of this system was inspired by 
Hooper~\cite{Hooper} and Schwartz~\cite{Schwartz}, but the specifics of our 
approach are different. We begin with an elementary
introduction to the theorems of Pappus and Steiner. The results of the
paper are described in Section~\ref{section.results} (on page~\pageref{section.results}) 
after the required notation is available. 

\subsection{Pappus's Theorem} \label{section.pappus.theorem} 
\begin{figure}
\includegraphics[width=10cm]{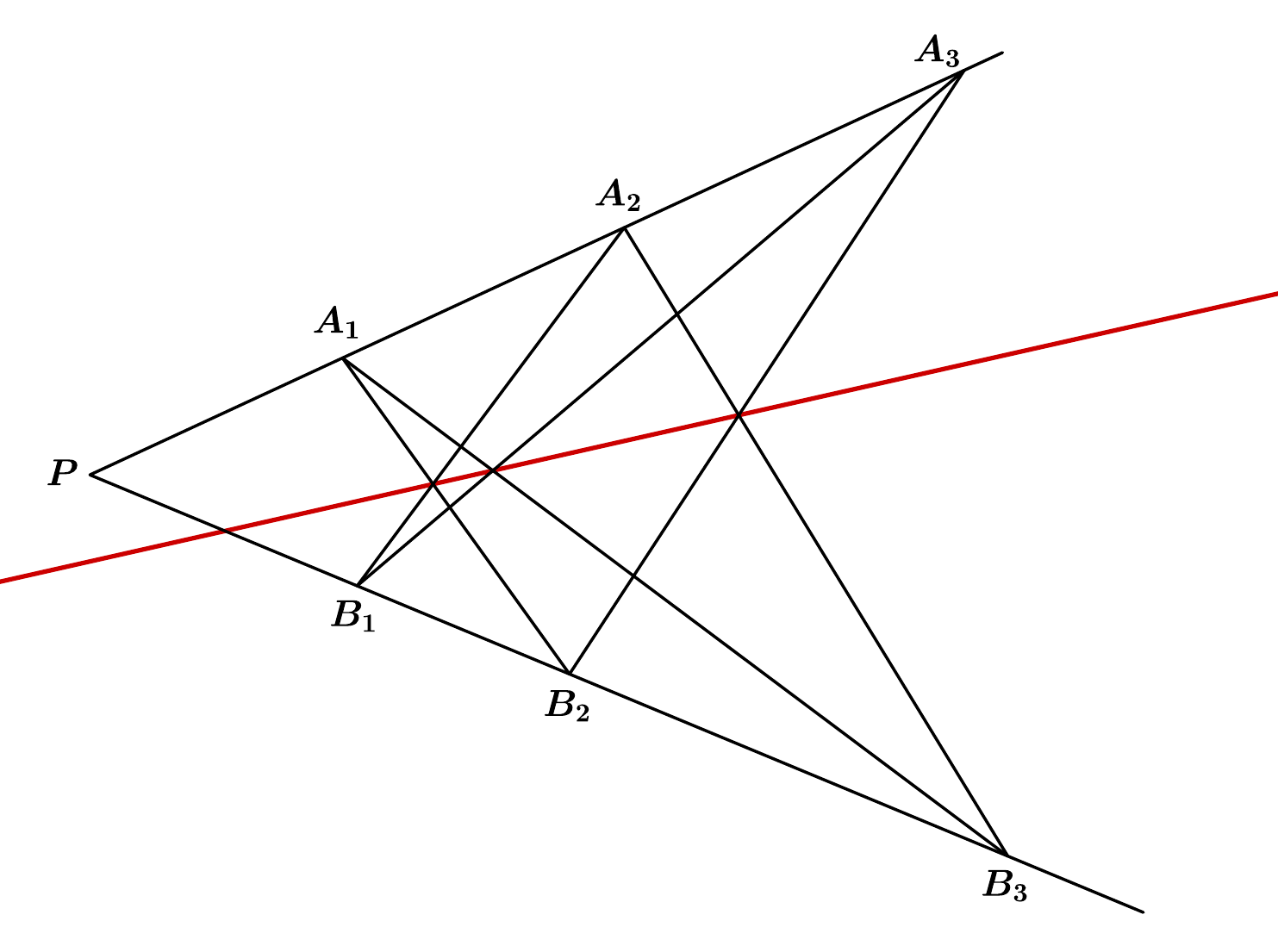} 
\caption{Pappus's theorem} 
\label{diag:pappus} 
\end{figure}

Let $\proj^2$ denote the projective plane over a
field\footnote{The field is allowed to be arbitrary for the moment, but we will make
  specific assumptions later.} $k$. Consider two lines 
$L$ and $M$ in $\proj^2$, intersecting at a point $P$. Choose distinct
points $A_1, A_2,A_3$ and $B_1, B_2, B_3$ on $L$ and $M$ respectively
(all away from $P$),  displayed as an array 
$\left[ \begin{array}{ccc} A_1 & A_2 & A_3 \\ B_1 & B_2 & B_3 \end{array} \right]$. Then Pappus's 
theorem says that the three cross-hair intersection points 
\[ A_1B_2 \cap A_2B_1, \quad A_2B_3 \cap A_3B_2, \quad A_1B_3 \cap A_3B_1, \] 
corresponding to the three minors of the array, are collinear (see
Diagram~\ref{diag:pappus}). The line containing them (called 
the Pappus line of the array) will be denoted by 
$\left\{\begin{array}{ccc} A_1 & A_2 & A_3 
\\ B_1 & B_2 &  B_3 \end{array} \right\}$. 

A proof of Pappus's theorem may be found in almost any book on
elementary projective geometry (for instance, see~\cite[Ch.~1]{Seidenberg}). 
\subsection{Steiner's Theorem} 
\label{section.steiner_theorem} 
\begin{figure}
\includegraphics[width=10cm]{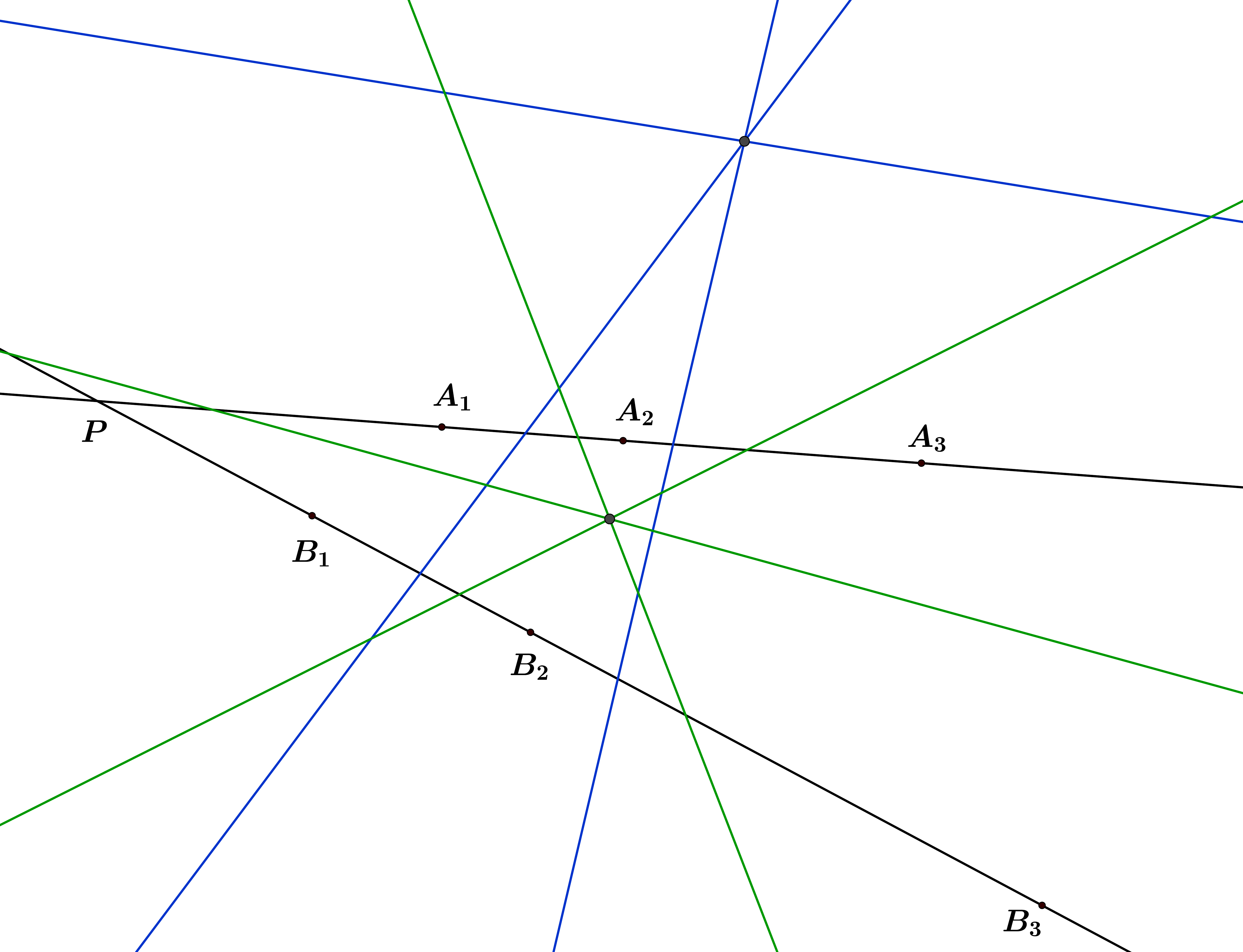} 
\caption{Steiner's theorem} 
\label{diag:steiner} 
\end{figure}

If we permute the bottom row of the array, then 
\emph{a priori} we get a different Pappus line. There are six such lines corresponding to the elements 
of the permutation group $\Sym_3$. For an element $\sigma \in \Sym_3$, let 
\[ \Lambda_\sigma = 
\left\{\begin{array}{ccc} A_1 & A_2 & A_3 \\ 
B_{\sigma(1)} & B_{\sigma(2)} & B_{\sigma(3)} \end{array} \right\} \] 
denote the corresponding Pappus line. 
According to our convention, the permutation $(1 \, 3 \, 2)$ takes $1$ to
$3$ etc., so that $\Lambda_{(1 \, 3 \, 2)} = \left\{\begin{array}{ccc} 
A_1 & A_2 & A_3 \\ B_3 & B_1 &  B_2 \end{array}  \right\}$. 
Now Steiner's theorem says that the 
lines $\Lambda_e, \Lambda_{(1 \, 2 \, 3)}, \Lambda_{(1 \, 3 \, 2)}$, 
corresponding to the even permutations are concurrent, and the lines 
$\Lambda_{(1 \, 2)}, \Lambda_{(1 \, 3)}, \Lambda_{(2 \, 3)}$
corresponding to the odd permutations are also concurrent. In
Diagram~\ref{diag:steiner}, the green lines correspond to the even permutations and blue lines to the odd
ones. In order to reduce visual clutter, the cross-hair lines have not been shown. 

A proof of Steiner's theorem may be found in~\cite[p.~216]{Baker} or~\cite[Section 5]{Hooper}. In
summary, starting from $A_1, \dots, B_3$, the entire process can be
seen as: 
\[ \text{Two lines with three points on each} \quad \leadsto \quad 
\text{Two points with three lines on each}. \] 
The latter structure (seen in the dual projective plane) is the same as the former, 
hence we will get a dynamical system (i.e., a map from a set to
itself) at an appropriate level of abstraction. We begin by recalling some standard facts about 
cross-ratios (see~\cite[Ch.~1]{Seidenberg}). 

\subsection{Cross-ratio and the $j$-function} \label{section.cross.j} 
Given a sequence of four points $x_1, x_2, x_3, x_4$ on the projective line 
$\proj^1$, define their cross-ratio to be 
\[ \langle x_1, x_2, x_3, x_4 \rangle = 
\frac{(x_1-x_3)(x_2-x_4)}{(x_2-x_3)(x_1-x_4)}. \] 
If $r$ denotes this value, then permuting the $x_i$ changes the cross-ratio into one
of the six possible functions of $r$. In fact, if $x_4$ is left fixed then the six permutations of $x_1, x_2, x_3$ give each of these functions exactly once. If $\sigma$ denotes any of the permutations 
\[ e, \quad (1 \, 2), \quad (1 \, 3), \quad (2 \, 3), 
\quad (1 \, 2 \, 3), \quad (1 \, 3 \, 2), \] 
then $\langle x_{\sigma(1)}, x_{\sigma(2)}, x_{\sigma(3)},
x_4\rangle$ is respectively equal to 
\begin{equation} 
r, \quad 1/r, \quad r/(r-1), \quad 1-r, \quad (r-1)/r, \quad 1/(1-r). 
\label{cross-ratio.values} \end{equation} 
Now define 
\[ j(r) = \frac{4 \, (r^2-r+1)^3}{27 \, r^2 \, (r-1)^2}.  \] 
This expression remains unchanged if we substitute any of the values in~(\ref{cross-ratio.values}) in place
of $r$. In other words, if $j(r) =
\tau$, then the six roots of the polynomial equation 
\[ 4 \, (Z^2-Z+1)^3 - 27 \, \tau \, Z^2 \, (Z-1)^2 = 0, \] 
are exactly those in (\ref{cross-ratio.values}). Henceforth we assume
that $\chr(k) \neq 2, 3$, since the theory of cross-ratios has
more pitfalls in those characteristics. 

\section{The Pappus structure} 
\subsection{} 
Define a Pappus structure in a projective plane to be an unordered pair of sets 
\begin{equation} 
\Pi = \left\{ \left\{A_1, A_2, A_3 \right\}, \left\{B_1, B_2, B_3
  \right\} \right\}, 
\label{pappus.structure1} \end{equation}
such that 
\begin{enumerate} 
\item the six points $A_1, \dots, B_3$ are pairwise
distinct, 
\item the three points in each set are collinear, 
\item the line containing the first three is different from the line containing
the other three, and 
\item the point of intersection of the two lines (say
$P$) is different from $A_1, \dots, B_3$. 
\end{enumerate} 

We will represent such a structure by a pair of numbers, called its
signature. The process followed in Steiner's theorem, when
applied to the signature, will give a dynamical system. 

\subsection{} Given a Pappus structure $\Pi$, consider the cross-ratios 
\[ r = \langle A_1, A_2, A_3, P\rangle, \qquad s = \langle B_1, B_2, B_3, P \rangle, \] 
and write 
\[ x = j(r)+j(s), \qquad y= j(r) \, j(s). \] 

\begin{Definition} \rm 
The ordered pair $[x,y]$ will be called the signature of $\Pi$,
denoted $\sig(\Pi)$. 
\end{Definition} 
It is well-defined by what we have said above. The intuition
behind this definition is two-fold: 
\begin{itemize} 
\item 
since the order of the three points within the line is
irrelevant, we pass from the cross-ratio to its $j$-function, and
\item 
since the two lines are on equal footing, we pass from $j(r)$ and $j(s)$ 
to their elementary symmetric functions. 
\end{itemize} 

Suppose that we have two Pappus structures $\Pi_1$ and $\Pi_2$, whose
points come from possibly different planes $\proj^2_{(1)}$ and
$\proj^2_{(2)}$. 
\begin{Definition} \rm 
We define $\Pi_1, \Pi_2$ to be equivalent, if there is an isomorphism 
$\proj^2_{(1)} \rightarrow \proj^2_{(2)}$ carrying $\Pi_1$ to
$\Pi_2$. 
\end{Definition} 
The following result justifies the definition of the signature. 

\begin{Proposition}\rm 
Two Pappus structures are equivalent, if and only if their signatures
coincide. 
\label{prop.equivalence} \end{Proposition} 
\proof The `only if' part is obvious, and the converse is a simple
computation (see Section~\ref{section.proof.equivalence}). \qed 

Let $\Pi$ be a Pappus structure with $\sig(\Pi) = [x,y]$. Assume $x \neq y$, and let 
\begin{equation} \steiner(\Pi) = \left\{ 
\left\{ \Lambda_e, \Lambda_{(1 \, 2 \, 3)}, \Lambda_{(1 \, 3 \, 2)}, \right\}, 
\left\{ \Lambda_{(1 \, 2)}, \Lambda_{(1 \, 3)}, \Lambda_{(2 \, 3)}
\right\} \right\}. 
\label{pappus.structure2} \end{equation}
We will show later (see Section~\ref{section.signature.proof}) that this is a
Pappus structure. Its `points' come from the dual
projective plane $(\proj^2)^\vee$. It is well-defined, since the starting points are
permuted in all possible ways during the construction. The next theorem
gives its signature. 
\begin{Theorem} \rm With notation as above, 
\[ \sig(\steiner(\Pi)) = \left[ \frac{2y \, (y-x+2)}{(x-y)^2},  \frac{y^2}{(x-y)^2} \right]. \] 
\label{theorem.signature} \end{Theorem} 
\proof See Section~\ref{section.signature.proof}. \qed 

Now consider the Pappus-Steiner map 
\begin{equation} 
\pi: k^2 \; - \rightarrow k^2, \quad 
[x,y] \longrightarrow \left[ \frac{2y \, (y-x+2)}{(x-y)^2},
\frac{y^2}{(x-y)^2} \right]. 
\label{formula.pappus_steiner} \end{equation}
Of course, $\pi$ is defined only if $x \neq y$. However, one can
define a subset $\affine \subset k^2$ such that $\pi$ and all of its iterates
are defined over it (see Section~\ref{pi.domain}). This leads\footnote{The approaches in~\cite{Hooper}
  and~\cite{Schwartz} retain the planar embedding of the
  configuration, whereas by contrast we distinguish it only  
up to projectivity.}  to a dynamical system $\pi: \affine \longrightarrow \affine$. 

\subsection{Results} \label{section.results} 
The system $\pi$ can be studied over any ground field, and either from
an analytic or an algebraic viewpoint. In this paper, we mostly do the
latter. 
\begin{enumerate}[labelindent=0pt] 
\item 
The basic geometric properties of $\pi$ are given in
Section~\ref{section.geometry.pi}. In particular, we show that the map
interchanges two geometrically natural subsets of $\affine$, 
namely the 'harmonic' and the 'balanced' Pappus structures. Moreover,
$\pi$ has a natural involution associated to it, which carries an 
interesting interpretation in terms of cross-ratios. 
\item 
One of the natural objects of study of a dynamical system is its
periodic points. In Theorems~\ref{theorem.period3} and
\ref{theorem.period4}, we characterise all (sufficiently large) primes $p$, such that $\pi$
admits periodic points of orders $3$ or $4$ over the prime field
$\field_p$. The technique involves a Gr{\"o}bner basis computation
followed by the use of class field theory. The polynomials which
emerge during the Gr{\"o}bner computation seem to be unusual from a
Galois-theoretic viewpoint, and this observation leads to a conjecture about the field
extension generated by all $n$-periodic points. 
\item 
The map $\pi$ has two fixed points and a $2$-cycle, and hence it is 
geometrically natural to consider the pencil of conics through these
four points (see Section~\ref{section.two_conics}). 
We discover that there is a \emph{unique} conic
in this pencil which is sent by $\pi$ into another such conic. This
leads to a one-dimensional quadratic dynamical system. 
In Section~\ref{section.complex_dynamics} we deduce a result about its
Julia set. 
\item 
Section~\ref{section.leisenring} contains a discussion of a theorem
due to Leisenring which is formally similar to Steiner's. We show that
it leads to exactly the same dynamical system as
in~(\ref{formula.pappus_steiner}). 
\item 
The appendix by Attila D{\'e}nes contains a brief computer-aided analysis of the
Pappus-Steiner system over real numbers. 
\end{enumerate} 
We will use~\cite{Coxeter, KK, Seidenberg} as standard references for
projective geometry. The reader is referred to~\cite{Lang, Milne} for the
necessary concepts from algebraic number theory. The basic terminology 
of dynamical systems may be found in~\cite{Devaney}, and that of Gr{\"o}bner bases
in~\cite{AL}. 

All the algebraic computations, including those
for Gr{\"o}bner bases, were done in~{\sc Maple} and confirmed
in~{\sc Macaulay-2}. The number-theoretic computations were done using a
combination of {\sc Maple} and Sage. For ease of reading, we have 
relegated some of the more computational proofs to
Section~\ref{section.proofs}.

\section{The geometry of the Pappus-Steiner map} 
\label{section.geometry.pi} 

\subsection{} \label{pi.domain} 
Let $W_1 = \{ [x,y]: x=y \}$ denote the locus in $k^2$ over which $\pi$ is not
defined. Then $\pi^2$ is not defined over the locus $W_2$ of points $[x,y]
\in k^2 \setminus W_1$ such that $\pi([x,y]) \in W_1$. In general,
define inductively 
\[ W_n = \{ [x,y] \in k^2 \setminus \bigcup\limits_{i=1}^{n-1} W_i: 
\pi([x,y]) \in W_{n-1} \} \quad \text{for $n \geqslant 2$},  \] 
and write $\affine = k^2 \setminus \bigcup\limits_{i=1}^\infty
W_i$. Then $\pi$ and all of its iterates are defined over $\affine$,
and we have a dynamical system $\pi: \affine \longrightarrow
\affine$. 

If $k$ is either the field of real or complex numbers, then $\affine$ is a   
dense open subset of $k^2$. 
A calculation shows that $W_2$ is the union of lines $y=0, y=2x-4$,
and $W_3$ is the curve $4x^2-4xy+y^2-8y=0$. There
seems to be no easy general formula for the equation of $W_n$.

\subsection{Balanced and Harmonic structures} 
\label{section.BH} 
We will say that a Pappus structure is balanced, if
$j(r)=j(s)$; which is equivalent to $x^2=4y$. This is tantamount to 
requiring that the unordered quadruples 
\[ \{P, A_1, A_2, A_3\}, \quad \{P, B_1, B_2, B_3\}, \] 
should be projectively isomorphic. Let $B = \{[x,y] \in \affine:
x^2=4y\}$ denote the parabola\footnote{It is understood that
  some points of the `parabola' are missing, since $\affine$ is a proper
  subset of $k^2$. The same will be true of the `line' $H$.} formed by balanced structures. 

Recall that four points on $\proj^1$ are said to be harmonic if their
cross-ratio (in some order) is $-1$. We will say that $\Pi$ is
harmonic, if either of the quadruples above is harmonic. This is
equivalent to either $j(r)$ or $j(s)$ being $=1$, i.e., $x=y+1$. 
Let $H = \{[x,y] \in \affine: x = y+1 \}$ denote the line formed by harmonic
structures. There is a simple but pleasing relation between these two
loci: 

\begin{Proposition} \rm 
With notation as above, $\pi(H) \subseteq B$ and $\pi(B) \subseteq H$. 
\label{prop.harmonic_balanced} \end{Proposition} 
\proof 
We have $\pi([t+1,t]) = [2t,t^2]$, which is balanced. Similarly, 
\[ 
\pi([2t,t^2]) = \left[\frac{2(t^2-2t+2)}{(t-2)^2}, \frac{t^2}{(t-2)^2}\right], \] 
which is harmonic. \qed 

Hence the second iterate $\pi^2$
sends $H$ to $H$ and $B$ to $B$. Both maps are given by the same
formula 
\begin{equation} \alpha(t) = \frac{t^2}{(t-2)^2}. 
\label{map.BH} 
\end{equation}
In Section~\ref{section.complex_dynamics}, we will treat this as a map
of one complex variable and analyse its dynamics. 

\subsection{} In order to investigate the image of $\pi$, we
write $\pi([x,y]) = [p,q]$, and try to solve for $x,y$ in terms of
$p,q$. This leads to equations 
\[ 2y(y-x+2)^2=p(x-y)^2, \quad y^2=q(x-y)^2. \] 
After eliminating $x$, we have 
\[ (p^2-4q) \, y^2-8pq \, y+16q^2 =0, \quad x =
\frac{(2q-p)y}{2q}+2. \] 
If $ p^2 \neq 4q$, then we get two solutions $y = \frac{4q}{p \pm 2
  \sqrt{q}}$, either of which determines $x$ uniquely. If $p^2 = 4q$,
then we get a unique solution 
\[ x =  \frac{2q}{p} +1, \quad y = \frac{2q}{p} \] 
in $H$. Hence we have the following proposition: 
\begin{Proposition} \rm 
The map $\pi: \affine \longrightarrow \affine$ is a double cover
ramified over $B$. 
\end{Proposition} 

This suggests that $\affine$ should have an \emph{involution}, i.e., a 
function $\affine \stackrel{\tau}{\longrightarrow} \affine$ which exchanges the
two sheets of the cover. In order to find it, let
$[x_1,y_1],[x_2,y_2]$ denote the two pre-images of $[p,q]$, so that 
we must have 
\[ x_i = \frac{(2q-p)y_i}{2q}+2 \quad (i=1,2), \quad 
y_1 + y_2 = \frac{8pq}{p^2-4q}, \quad y_1 \, y_2 =
\frac{16q^2}{p^2-4q}. \] After eliminating $p$ and $q$, we get 
\[ x_2 = \frac{x_1-2y_1}{x_1-y_1-1}, \quad y_2 =
\frac{y_1}{x_1-y_1-1}. \] 

Hence we have the formula 
\[ \tau([x,y]) = \left[\frac{2y-x}{y-x+1},\frac{y}{y-x+1}\right]. \] 
As it stands, it is only defined over $\affine \setminus H$. However, it is natural to set 
$\tau([x,y]) = [x,y]$ for a point in $H$. 

\subsection{} \label{section.circ} 
This involution can be interpreted at the level of Pappus
structures. To see this, define 
\[ w^\circ = \frac{2-w}{1+w} \] 
for $w \neq -1, 2, 1/2$. This operation is such that 
$(w^\circ)^\circ =w$.  The next proposition shows that it imitates the action of $\tau$ at the level of
cross-ratios. 
\begin{Proposition} \rm 
Let $[x,y] = [j(r) + j(s), j(r) \, j(s)]$ denote a point away from $H$. Then 
\[ \tau([x,y]) = \left[ j(r^\circ)+j(s^\circ), j(r^\circ) \, j(s^\circ) \right]. \] 
\end{Proposition} 
\proof
Consider the equation $Z^2 - \frac{2y-x}{y-x+1} \, Z + \frac{y}{y-x+1}
=0$. Applying the quadratic formula, one of its roots is 
\[ \frac{4 \, (r^2-r+1)^3}{(r-2)^2 \, (2r-1)^2 \, (r+1)^2}, \] 
and of course, the second root is the same expression in $s$. If
we write $j(z)$ for this expression, then it leads to a sextic
polynomial equation in $z$ with roots 
\begin{equation} 
\frac{2-r}{r+1}, \quad \frac{r+1}{2-r}, \quad \frac{r+1}{2r-1}, \quad \frac{2r-1}{r+1},
\quad \frac{2r-1}{r-2}, \quad \frac{r-2}{2r-1}. 
\label{tau.cr_values} \end{equation} 
We declare the first of these to be $r^\circ$, but the choice is
arbitrary. One can get them all by applying the functions in~(\ref{cross-ratio.values})
to $r^\circ$. \qed

The relationship between the $j$-functions of $r$ and $r^\circ$ is also
involutive. We have 
\[ j(r^\circ) = \frac{j(r)}{j(r)-1}, \quad \text{and} \quad j(r) =
\frac{j(r^\circ)}{j(r^\circ)-1}. \] 
This follows by a direct computation using the formula for
$r^\circ$. 
\section{Periodic points} \label{section.periodic_points} 

\subsection{} 
Let $n > 0$ be an integer. Recall that a point $z=[x,y] \in \affine$
is said to be $n$-periodic, if $\pi^n(z) = z$. The smallest integer $n$ for which this
equation holds is called the period of $z$. A point of period $1$ is
called a fixed point. If $z$ is $n$-periodic, then 
\[ z \rightarrow \pi(z) \rightarrow \pi^2(z) \rightarrow \dots
\rightarrow \pi^n(z)=z, \] 
is called an $n$-cycle. 

Since $\pi$ is expressed by rational functions, it is natural to use
Gr{\"o}bner bases to detect the points of period $n$. Fix the polynomial ring $R =
\QQ[x,y]$ with a lexicographic term order given by $x > y$. This is
suitable for eliminating the variable $y$ (see~\cite[Ch.~2]{AL}). 

For illustration, assume $n=2$. Then the equation $\pi^2(z) = z$ gives the ideal 
\[ I_2 = (y(4x^2-4xy+y^2-16x+7y+16), 
-4x^3+4x^2y-xy^2+20x^2-12xy+2y^2-16x-8y) \] 
in $R$. A point $[x,y]$ is $2$-periodic, if and only if both polynomials in
$I_2$ vanish on it. The Gr{\"o}bner basis of this ideal is given by 
\[ \begin{aligned} 
f_1(x) & = x(x-1)(x-2)(x-4)(x-5)(x-8)(x-12), \\ 
f_2(x,y) & = 18480y-67x^6+1934x^5-19615x^4+84620x^3-152728x^2+85856x. 
\end{aligned} \] 
Now, although this computation was carried out over $\QQ$, this is
also the Gr{\"o}bner basis of $I_2$ (seen as an ideal in $k[x,y]$) 
for $\chr(k)$ sufficiently large. This gives the following proposition: 

\begin{Proposition} \rm 
Assume that $\chr(k)$ is either zero or sufficiently large. Then 
the only fixed points of $\pi$ are $[2,1]$ and $[12,16]$. Moreover, 
the only points of period $2$ are $[5,4]$ and $[8,16]$, which are sent
to each other by $\pi$. 
\end{Proposition} 
\proof This follows by solving $f_1$ for $x$, and then
using $f_2$ to find $y$. The values $x=0,1$ lead to $y=0$, which is
disallowed. \qed 

We are aiming for a similar theorem for the cases $n=3, 4$. Assume
that $k = \field_p$, for $p$ sufficiently large. 
\subsection{Case $n=3$} \label{section.period3} 
We can find the Gr{\"o}bner basis of $I_3$ as before. It turns out to
be 
\[ f_1(x) = x^9(x-12)(x-2)
\underbrace{(x^3-84x^2+896x-1984)}_{u(x)} 
\underbrace{(x^3-48x^2+656x-1856)}_{v(x)}, \] 
(where $u(x), v(x)$ are irreducible over $\QQ$) together with 
\[ f_2(x,y) = x^4 \, y-\phi_1(x), \qquad f_3(x,y) = y^2+yx-\phi_2(x). \] 
Here $\phi_1$ and $\phi_2$ are degree $16$ polynomials in $x$,
which are too cumbersome to write down in full. Then one checks that 
$x=0$ is not an admissible value, and $x=2$ and $12$ lead to the fixed
points which we have already found. Moreover, if we substitute the value
of $y$ obtained from $f_2$ into $f_3$, then the latter vanishes
identically modulo $f_1$. Hence there exists a point of period $3$, 
exactly when either $u(x)$ or $v(x)$  has a root in $\field_p$.  

Now the pleasant surprise is that the discriminants of these
polynomials, namely 
\[ \Delta_u = 2^{12} \, 13^2 \, 31^2 \quad \text{and} \quad 
\Delta_v = 2^{12} \, 7^2, \] 
are both complete squares, which implies that their Galois groups over
$\QQ$ are both isomorphic to the cyclic group $C_3$. Hence we can apply class field
theory to find conditions on $p$ such that either of these polynomials admits a root modulo
$p$. The details will be given in
Section~\ref{section.class_field.theory}, but here we state the resulting theorem: 

\begin{Theorem} \rm The system $\pi$ has no points of period $3$ over
  $\QQ$. For all sufficiently large primes $p$, it has such a point
  over $\field_p$ if and only if, 
\[ p \equiv 1, 5, 8, 12 \mod{13}, \quad \text{or} \quad  
p \equiv 1, 6  \mod{7}.  \] 
\label{theorem.period3} \end{Theorem} 
For example, if $p =47$, then 
$[21,2] \rightarrow [39,6] \rightarrow [24,28] \rightarrow [21,2]$ is a $3$-cycle. 

\subsection{Case $n=4$} The argument is similar. The Gr{\"o}bner basis
of $I_4$ consists of two polynomials: 
\[ \begin{aligned} 
f_1(x)  =  & \; x \,  (x-1) \, (x-2) \, (x-4) \, (x-5) \, (x-8) \,
(x-12) \times \\
&  (x^2-96x+304) \,  (x^2-14x+29) \, (x^2-20x+80) 
(x^2-12x+16) \, (x^2-24x+64) \times \\
& \underbrace{(x^4-144x^3+3152x^2-16896x+25856)}_{w(x)}, \end{aligned}
\] 
together with $f_2(x,y) = y - \psi(x)$, where 
$\psi(x)$ is a degree $20$ polynomial which need not be
written down explicitly. The linear factors in $f_1(x)$ lead to the fixed points or $2$-cycles
already found. Hence, as before, we have a point of period $4$ in
$\field_p$ if and only if, either any of
the quadratic factors or the unique quartic factor has a root. 

The quadratics are easy to take care of. For instance, 
$x^2-96x+304$ has a root in $\field_p$ exactly when its discriminant
$2^6 \, 5^3$ is a square in $\field_p$, i.e., if and only if the
Legendre symbol $\left(\frac{5}{p}\right)=1$. By quadratic
reciprocity, the last condition is equivalent to $p \equiv 1, 4
\mod{5}$. The discriminants of the remaining quadratics are all of the
form `$\text{integer squared} \times 5$', hence they all lead to the same
criterion. 

Now, just as in the previous case it turns out that the quartic $w(x)$ has Galois group $C_4$,
hence we can use class field theory. The details will be given in
Section~\ref{section.cases.3and4}; here we state the result: 
\begin{Theorem} \rm 
The system $\pi$ has no points of period $4$ over
  $\QQ$. For all sufficiently large primes $p$, it has such a point
  over $\field_p$ if and only if, 
\[ p \equiv 1, 4 \mod{5}, \quad \text{or} \quad  
p \equiv 1, 4, 13, 16  \mod{17}.  \] 
\label{theorem.period4} \end{Theorem} 
For example, if $p=89$, then 
$[80,36] \rightarrow [8,22] \rightarrow [49,17] \rightarrow [7,44]
\rightarrow [80,36]$ is a $4$-cycle. 

\subsection{A Galois-theoretic conjecture} 
Let us recapitulate the case $n=4$. The polynomial $f_1(x)$
has several irreducible factors over $\QQ$, all of which have abelian
Galois groups. Hence the compositum $K$ of their splitting fields is also an abelian extension
of $\QQ$. In fact, if $\theta$ denotes any of the roots of $w(x)$, then $K =
\QQ(\theta, \sqrt{5})$, which has Galois group $C_4 \times
C_2$. Moreover, the structure of $f_2$ is such that $y$ is completely determined as a rational function of $x$. 

It is a very surprising fact that the corresponding statements are also true of
$n=3,5,6$. For instance,  when $n=5$, the polynomial $f_1(x)$ has
one irreducible factor of degree $5, 10$ and $15$ each (not counting
the linear factors), and all three have abelian Galois groups.\footnote{See~\cite{Hulpke} for a survey of
  algorithms for computing Galois groups of polynomials over
  $\QQ$. In practice, computer algebra systems tend to use 
opportunistic combinations of such techniques.}  
Since this is extremely unlikely for a `random' choice of polynomials
of such high degrees,
it is all but certain that the phenomenon must be a manifestation of a deeper
hidden structure. We present this as a conjecture. 

\begin{Conjecture} \rm 
Assume $n >2$, and let $\{[\alpha_i, \beta_i]: i =1, \dots, r\}$ be the set of
$n$-periodic points of $\pi$ over $\QQ$. Let 
$K = \QQ(\{\alpha_i, \beta_i\})$ denote the number field generated by
all of their coordinates. Then $K/\QQ$ is an abelian extension. 
\label{conjecture.abelian} \end{Conjecture} 

The Galois group of $K/\QQ$ is 
\[ C_3 \times C_3, \quad C_2 \times C_4, \quad C_2 \times C_3 \times C_5^3, \] 
for $n=3,4,5$ respectively. The case $n=6$ already seems too large for
a complete answer, but a partial computation shows that it must a
subgroup of $C_3^8 \times C_4$. 

If the conjecture is true, then it may well be possible to 
find similar congruence conditions for all $n$. But this seems out of reach at the moment. 

\section{The two conics} 
\label{section.two_conics} 

The Pappus-Steiner map turns out to have an intriguing and unexpected
relationship with a pair of conics in the projective plane. We begin
by explaining the underlying geometric construction, which culminates in
Theorem~\ref{theorem.2conics}. Throughout this section, we work
over the complex numbers. 

\subsection{} Let $\Gamma$ denote the projective plane over $\complex$, with
homogeneous coordinates $x,y,z$. Then $\affine$ can be seen as a subset
of $\Gamma$ via the natural inclusion\footnote{Although
  $\Gamma$ and $\proj^2$ are isomorphic, they play conceptually
  different roles. In a sense, $\Gamma$ is an `extension' of
  $\affine$, whereas $\proj^2$ is the natural home for configurations
  which are parametrised by $\affine$. In this section 
  $\proj^2$ will not play any role.} 
\[ \affine \longrightarrow \Gamma, \quad [x,y] \longrightarrow [x,y,1]. \] 
Now consider the homogenized version of $\pi$ (denoted by the same
letter) 
\[ \pi: \Gamma \longrightarrow \Gamma, \quad [x,y,z] \longrightarrow 
\left[ 2y(y-x+2z), y^2, (x-y)^2 \right];  \] 
together with the involution 
\[ \tau: \Gamma \longrightarrow \Gamma, \quad [x,y,z] \longrightarrow
[2y-x,y,y-x+z]. \] 
Of course, $\pi$ and $\tau$ are only rational maps, i.e., they are
defined over a dense open subset of $\Gamma$. 
\subsection{} 
Now consider the points
\[ R_1 = [2,1,1], \quad R_2 = [12,16,1], \quad R_3 = [5,4,1], \quad
R_4 = [8,16,1] \] 
in $\Gamma$. Recall that $R_1, R_2$ are fixed points of
$\pi$ and $R_3 \leftrightarrow R_4$ form a $2$-cycle. Since $\pi$
stabilises the set $\{R_1, \dots, R_4\}$, it is natural to study its 
behaviour with respect to the pencil of conics passing through the
$R_i$. 

One can describe this pencil as follows (cf.~\cite[Ch.~6]{Seidenberg}). 
If $\conic$ is any smooth conic through $R_1,
\dots, R_4$, then as a point $P$ moves on $\conic$, the cross-ratio of the four
lines $PR_1, PR_2, PR_3, PR_4$ remains independent of $P$. Moreover, this
common value identifies the conic. Thus, for any
$q$, we define\footnote{This description breaks down for the three
  degenerate conics in the pencil. Moreover, when $P$ approaches any
  of the $R_i$, one must take the limiting value of the
  cross-ratio. These subtleties will cause no problems for
  us.} the conic 
\[ \conic_q = \{P \in \Gamma: \langle PR_1, PR_2, PR_3, PR_4 \rangle
=q \},  \] 
passing through $R_1, \dots, R_4$. 

For a parameter $t$, consider the line $L_t = \{[x,y,z]: tx-y+(1-2t)z=0 \}$
passing through $R_1$. It is natural to set $L_\infty$ to be the line
$x-2z=0$. Let us identify $\proj^1$ with the pencil of lines through $R_1$ via the map 
\[ t \longrightarrow L_t. \] 
Then we have an isomorphism 
\[ \proj^1 \longrightarrow \conic_q, \quad t
  \longrightarrow P_q(t), \] 
where the conic $\conic_q$ and the line $L_t$ intersect in the point-pair $R_1$ and $P_q(t)$. 

\subsection{} 
It is clear that $\pi(\conic_q)$ must be an algebraic curve passing
through $R_1, \dots, R_4$, but unfortunately it is not a conic in
general. However, for the~\emph{unique} value $q =2/7$,
it is in fact the conic $C_{-1/4}$ belonging to the same pencil. This
situation gives a pretty geometric example. In order to state the
result in full, define the two functions 
\[ \imath(t) = \frac{19t-23}{16t-19}, \quad \text{and} \quad 
\beta(t) = \frac{44t^2-76t+27}{56t^2-120t+62}. \] 
Notice that $\imath$ is involutive, i.e., $\imath^2(t) =t$. Since
$\imath$ is a M{\"o}bius transformation, this corresponds to the
identity $\left[\begin{array}{rr} 19 & -23 \\ 16 & -19 \end{array}
\right]^2 = -7 \, \left[\begin{array}{rr} 1 & 0 \\ 0 &
                                                       1 \end{array}\right]$. 

\begin{figure}
\includegraphics[width=12cm]{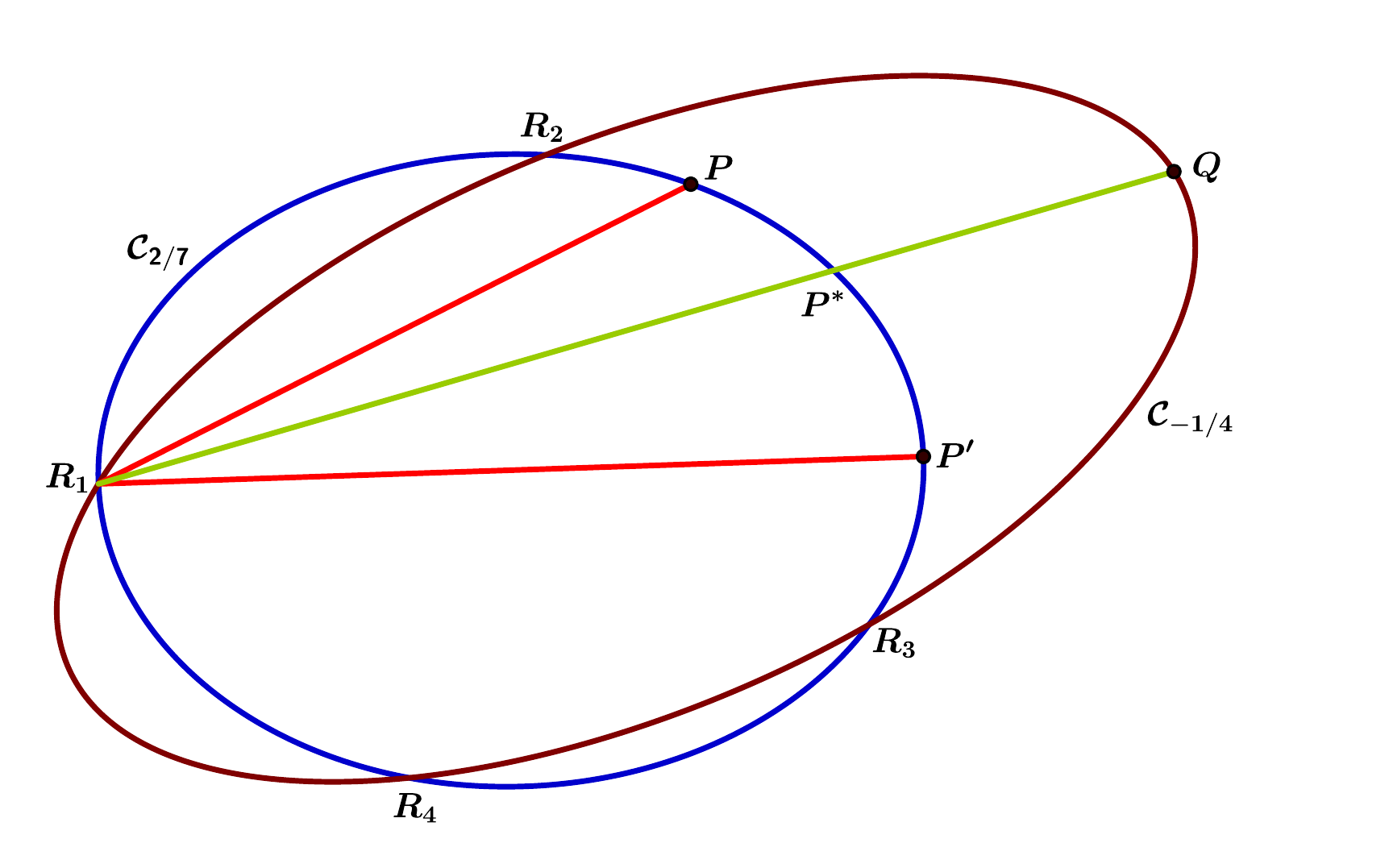} 
\caption{The double cover $C_{2/7} \stackrel{\pi}{\longrightarrow} C_{-1/4}$} 
\label{diag:conics} 
\end{figure}

\begin{Theorem} \rm 
We have 
$\tau(\conic_{2/7}) = \conic_{2/7}$ and $\pi(\conic_{2/7}) =
\conic_{-1/4}$. That is to say, the conic $\conic_{2/7}$ is stabilised
as a set by the involution $\tau$, and mapped onto $\conic_{-1/4}$ by $\pi$. Moreover, 
\[ \tau(P_{2/7}(t)) = P_{2/7}(\imath(t)), \quad \text{and} \quad 
\pi(P_{2/7}(t)) = P_{-1/4}(\beta(t)). \] 
\label{theorem.2conics} \end{Theorem} 
The proof follows by a computation, which is given in
Section~\ref{section.conic_computation}. The geometry behind the
theorem can be seen in Diagram~\ref{diag:conics}. 

\begin{itemize} 
\item 
The blue conic $\conic_{2/7}$, and the brown conic $\conic_{-1/4}$ 
both pass through $R_1, \dots, R_4$. 
\item 
Given any point $P$ on $\conic_{2/7}$, the involution $\tau$ sends it to
a point $P'$ on the same conic. Moreover, $\pi(P) = \pi(P') = Q$ always lies on
$\conic_{-1/4}$. 
\item 
If points on either conic are described by the parameter $t$ as above, then the
function $P \rightarrow P'$ corresponds to $t \rightarrow \imath(t)$ and $P
\rightarrow Q$ corresponds to $t \rightarrow \beta(t)$. In short, $\imath$ imitates
the action of $\tau$ and $\beta$ imitates the action of $\pi$. 
\end{itemize} 

The set-up of Steiner's theorem contains no hint of such an example, hence its
existence is something of a curiosity. Although the two conics in
the theorem are intrinsically determined, the parameter values $2/7$ and $-1/4$ depend on the
choice of $R_1$ in defining the point $P_q(t)$. 
\subsection{Complex Dynamics} 
\label{section.complex_dynamics} 
In the next two sections we will analyse the complex dynamics of the
map $\beta(t)$, as well as that of $\alpha(t)$ from Section~\ref{section.BH}. With $\proj^1$
seen as the Riemann sphere, $\alpha(t)$ and $\beta(t)$ are quadratic rational maps 
$\proj^1 \rightarrow \proj^1$. 

An excellent introduction to one-dimensional complex dynamics 
may be found in~\cite{Beardon}. A detailed discussion of the
dynamics of quadratic rational maps is given in~\cite{Milnor}. 

Recall that $\alpha(t) = \frac{t^2}{(t-2)^2}$, and it gives the
formula for either of the maps $H \stackrel{\pi^2}{\longrightarrow} H, B
\stackrel{\pi^2}{\longrightarrow} B$. Now $\alpha$ has fixed points $0,1, 4$,
where the first is super-attracting and the other two are repelling. After the change of variables
$t = 2/(1+w)$ (which amounts to conjugating by a M{\"o}bius transformation),
the function takes the form $\widetilde{\alpha}(w) = 2w^2-1$. The latter
is the Tchebychev polynomial $T_2$, whose Julia set is known to be
the real closed interval $[-1,1]$ (see~\cite[\S 1.4]{Beardon}). 
Reverting to the $t$-variable, we deduce that the Julia set
$J(\alpha)$ is the real closed interval $[1,\infty]$. 

The basin of attraction of the fixed point $0$ is $\proj^1 \setminus
[1,\infty]$. In geometric terms, this means that if we start with a harmonic structure $[t+1,t]$
where $t$ lies in this basin, then for large $n$ the iterate 
$\pi^{2n}([t+1,t])$ will tend to a Pappus structure where one of the
lines has a point pair coming together. On the other hand, the 
trajectory is chaotic if $t \in [1,\infty]$. The same is true of
the balanced structure $[2t,t^2]$. 

\subsection{} The map $\beta(t)$ has fixed points 
\[ a = \frac{10-\sqrt{37}}{14} \approx 0.2798, \quad 
b = \frac{10+\sqrt{37}}{14} \approx 1.1488, \quad 
c = 3/2. \] 
The first is attracting, and the latter two are repelling. There is
also a repelling $2$-cycle $1 \leftrightarrow 5/2$. There are
two complex conjugate critical points $\frac{19 \pm i
  \sqrt{7}}{16}$. Now it follows from~\cite[Lemma 10.2]{Milnor} that the
Julia set $J(\beta)$ must be a Cantor set contained in $\RR$. In
principle, it would be possible to give a recipe for
constructing it just as for the usual 'middle-third' Cantor set 
(cf.~\cite[\S 1.8]{Beardon}), but we
omit this for brevity. 

The attracting fixed point has the following geometric
interpretation. If we start with a point $P$ as in
Diagram~\ref{diag:conics}, replace it with $P^*$ and 
continue this process indefinitely, then it 
will converge to $P_{2/7}(a) \approx [6.5951, 2.2857,1]$ for
almost all starting positions. 

\section{Leisenring's theorem} \label{section.leisenring} 
We now give a short account of a theorem
due\footnote{The question of priority seems a little muddled
  (see~\cite[\S 1]{Rigby}), since apparently Leisenring never
  published his result.}  to
Leisenring which is formally similar to Pappus's. It has an extension due to Rigby, which is similar to Steiner's. By the same considerations as in
Section~\ref{section.steiner_theorem}, we get a dynamical
system. But it turns out to be the same as the Pappus-Steiner
system, and as such does not lead to anything new. 

\subsection{} 
\begin{figure}
\includegraphics[width=10cm]{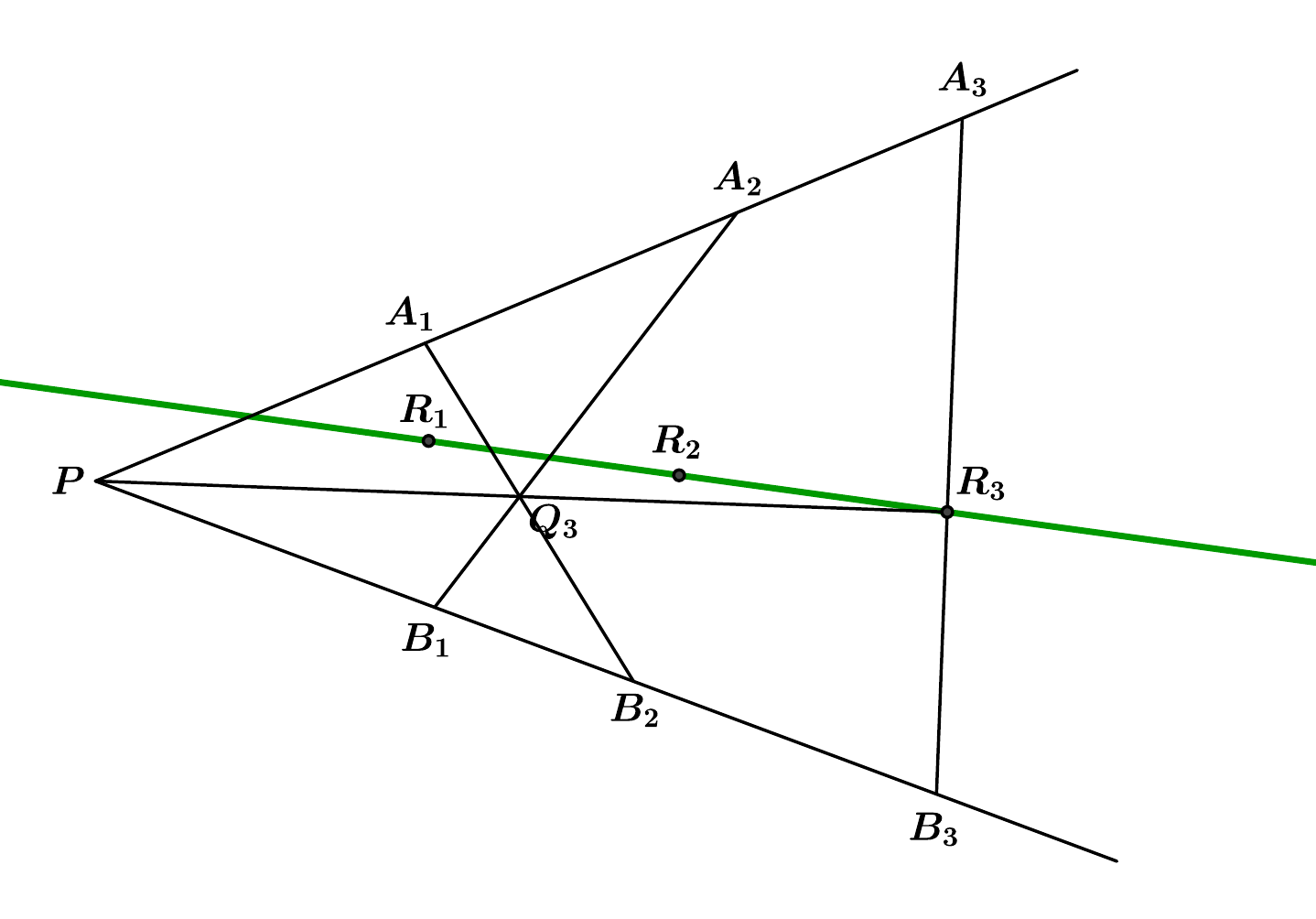} 
\caption{Leisenring's theorem} 
\label{diag:leisenring} 
\end{figure}

In the notation of Section~\ref{section.pappus.theorem}, define points 
\[ Q_1 = A_2B_3 \cap A_3B_2, \quad R_1 = PQ_1 \cap A_1B_1, \] 
and similarly\footnote{Of course, these points have nothing to do
  with the $R_1, \dots, R_4$ from the previous section.} for $R_2, R_3$. Then Leisenring's theorem says that 
$R_1, R_2, R_3$ are collinear (see Diagram \ref{diag:leisenring}). 
If $\left\langle\begin{array}{ccc} 
A_1 & A_2 & A_3 \\ B_1 & B_2 & B_3 \end{array} \right\rangle$ denotes
the Leisenring line containing the $R_i$, then we have the further
theorem due to Rigby that the six lines 
\[ \Psi_\sigma = \left\langle\begin{array}{ccc} 
A_1 & A_2 & A_3 \\ 
B_{\sigma(1)} & B_{\sigma(2)} & B_{\sigma(3)} \end{array}
\right\rangle, \quad \sigma \in \Sym_3 \] 
are also concurrent in threes exactly as in Steiner's theorem. The proofs
may be found in~\cite{Rigby}. We omit the 
diagram, since as a visual pattern it is no different from the one for 
Steiner's theorem. Let 
\[ \leisen(\Pi) = \left\{ 
\left\{ \Psi_e, \Psi_{(1 \, 2 \, 3)},  \Psi_{(1 \, 3 \, 2)} \right\}, 
\left\{ \Psi_{(1 \, 2)}, \Psi_{(2 \, 3)},  \Psi_{(1 \, 3)}
\right\} \right\}. \] 

\begin{Theorem} \rm 
Let $\sig(\Pi) = [x,y]$. If $x \neq y$, then $\leisen(\Pi)$ is a Pappus structure. Moreover, it is
equivalent to the structure $\steiner(\Pi)$ obtained from Steiner's theorem. 
\label{theorem.leisenring_steiner} \end{Theorem} 
Essentially the same theorem is proved in \cite[\S 3]{Rigby} by
synthetic methods. We will give an analytic proof in Section~\ref{section.leisen.proof}. In general, 
\[ [\Psi_\sigma: \sigma \in \Sym_3] \quad \text{and} \quad 
[\Lambda_\sigma: \sigma \in \Sym_3], \] 
are entirely different lines, but there is an automorphism
of $\proj^2$ carrying the first sequence to the second. We will find an
explicit matrix for it immediately after the proof. 

\section{Some computations} \label{section.proofs} 

We have collected some of the more computational proofs in
this section. 
\subsection{The `standard' Pappus structure} \label{section.constructionAB}
For any constants $r$ and $s$ (different from $0,1$) construct a collection of points
$C(r,s) = \{A_1, \dots, B_3 \}$ as follows. Choose coordinates $z_1, z_2, z_3$ in $\proj^2$, 
and let $L,M$ denote the lines $z_1=0, z_2=0$ respectively. Then 
$P=[0,0,1]$ is their point of intersection. The points 
\begin{equation} \begin{array}{lll} 
A_1 = [0,1,r], & A_2 = [0,1,1], & A_3 = [0,1,0], \\ 
B_1 = [1,0,s], & B_2 = [1,0,1], & B_3 = [1,0,0] 
\end{array} \label{standard.points} \end{equation}
are such that $\langle A_1, A_2, A_3, P \rangle =r$, and $\langle B_1, B_2, B_3, P 
\rangle =s$. Let $\Pi(r,s)$ denote
the Pappus structure determined by $C(r,s)$. It will prove useful in many of the computations below. 

If $N$ is a nonsingular $3 \times 3$ matrix, define an automorphism $\varphi_N : \proj^2 \longrightarrow 
\proj^2$ by sending the row-vector $[z_1 \; z_2 \; z_3]$ to $[z_1 \;
z_2 \; z_3] \, N$. 
\subsection{Proof of Proposition~\ref{prop.equivalence}} \label{section.proof.equivalence}
Write $r \sim r'$, if $r'$ is any of the functions of $r$
in~(\ref{cross-ratio.values}). 

\begin{Lemma} \rm 
Assume that $r \sim r'$ and $s \sim s'$. Then $\Pi(r,s)$ is isomorphic 
to $\Pi(r',s')$. 
\label{lemma.pappus.equivalence} \end{Lemma} 
\proof 
We want to find an automorphism $\varphi_N$ of $\proj^2$ which takes
$\Pi(r,s)$ to $\Pi(r',s')$. Let $N = \left[ \begin{array}{ccc} a & 0 & b \\ 0 & c & d 
\\ 0 & 0 &  1 \end{array} \right]$ for some $a,b,c,d$. Now it is easy to see that 
$\varphi_N$ preserves the lines $L$ and $M$, and hence also the point $P$. The image of 
a point on $L$ (resp.~$M$) depends only on $c, d$ (resp.~$a,b$). The
gist of the lemma is that, depending on $r'$ and $s'$, we can independently 
specify the images of $A_1, A_2, B_1, B_2$ under $\varphi_N$ by
choosing $a, \dots, d$ correctly. Then the standard properties of the
cross-ratio will ensure that $A_3$ and $B_3$ automatically land where we want them to. 

For instance, let $r = 1/(1-r')$ and write $A_1' =[0,1,r']$. Now 
$\langle A_1', A_2, A_3, P \rangle = r'$ implies that 
$\langle A_3, A_1', A_2, P \rangle = r$. If we can arrange 
the entries in $N$ such that $\varphi_N$ takes $A_1, A_2$ respectively to $A_3, A_1'$, 
then it will necessarily take $A_3$ to $A_2$. But this amounts to 
choosing $c,d$ such that $d+r=0, 1+d=c \, r'$, which can surely 
be done. Parallely, choose $a, b$ depending on $s'$. 
All the other cases are very similar, and we leave them to the reader. \qed 

\medskip 

Given a Pappus structure, we can always choose coordinates in the 
plane such that it is of the form $\Pi(r,s)$ for some $r,s$. Moreover,
$\Pi(r_1,s_1)$ and $\Pi(r_2,s_2)$ have the same signatures, if and 
only if, either $r_1 \sim r_2, s_1 \sim s_2$ or $r_1 \sim s_2, s_1
\sim r_2$. This proves Proposition~\ref{prop.equivalence}. \qed 

\subsection{Proof of Theorem~\ref{theorem.signature}}
\label{section.signature.proof} 
All the computations below were done in {\sc Maple}, and it is a pleasant surprise that several
intermediate expressions turn out to have clean and tidy factorisations. For the
sake of brevity, we will give only the salient steps in the proof. 

Starting from the standard Pappus structure, it is straightforward to calculate the 
coordinates of all Pappus lines. Let $E$ (resp.~$E'$) denote the point
of concurrency of all even (resp.~odd) lines. Their coordinates
are 
\[ \begin{aligned} 
E & = [r (r-1) (s^2-s+1), s \, (s-1) \, (r^2-r+1), rs(rs-1)], \\ 
E' & = [r(r-1)(s^2-s+1), -s(s-1)(r^2-r+1), rs(r-s)]. 
\end{aligned} \] 
Either of these points is not well-defined, if all three lines passing through it coincide.  
We claim that at least one of these points is not defined,  if and only 
if 
\begin{equation} r^2-r+1=0 \quad \text{and} \quad s^2-s+1=0. 
\label{condition.rs} \end{equation} 
The `only if' part is clear. For the converse, if $r, s$ are the roots 
of $Z^2 - Z+1=0$, then either they are equal (when $E'$ is 
undefined) or they are reciprocals (when $E$ is undefined). 
Henceforth we can assume that condition~(\ref{condition.rs}) is not 
satisfied. Then the two points are defined and turn out to be always 
distinct. Let $\Lambda_*$ denote the line joining them, and calculate the 
cross-ratios 
\[ c = \langle \Lambda_e, \Lambda_{(1 \, 2 \, 3)}, \Lambda_{(1 \, 3 
  \, 2)}, \Lambda_* \rangle, \qquad 
c' = \langle \Lambda_{(1 2)}, \Lambda_{(1 \, 3)}, \Lambda_{(2 \,
  3)}, \Lambda_* \rangle. \] 
We will omit the expressions for $c$ and $c'$, but write down
their $j$-functions instead. Let 
\[ \begin{aligned} 
u & = 4 \, (r^2-r+1)^3 \, (s^2-s+1)^3, \\ 
v & = (rs-2r+s+1)(rs+r-2s+1)(2rs-r-s+2), \\ 
v' & = (rs+r+s-2)(2rs-r-s-1)(rs-2r-2s+1). \end{aligned} \] 
Then $j(c) = u/v^2$ and $j(c') = u/v'^2$. 
These are well-defined, if and only if $v, v'$ are both nonzero. Now 
we have 
\[ x - y = j(r)+ j(s) - j(r) \, j(s) = \frac{-4 \, v \, v'}{[27 \, r (r-1) s (s-1)]^2}, \] 
so this happens exactly when $x \neq y$. 

Let $\sig(\steiner(\Pi)) = [p,q]$, which leads to equations 
\[ j(r) + j(s)-x = 0, \quad j(r) \, j(s)-y=0, \quad 
j(c) + j(c')-p=0, \quad j(c) \,  j(c')-q=0. \] We should like to find
$p,q$ in terms of $x,y$. Now consider the polynomial ring in variables $x,y,p,q,r,s$, with an 
elimination order with respect to the last two variables, and  
find the Gr{\"o}bner basis of the ideal generated by the numerators 
of these equations. Then $r$ and $s$ get eliminated, leading 
to \[ 2y \, (y-x+2) - p \, (x-y)^2 =0, \quad  \text{and} \quad y^2 - q \, (x-y)^2=0. \] 
This gives the required formula for the Pappus-Steiner map. 
The expressions show, once again, that 
$p,q$ are defined only when $x \neq y$. \qed 

We can identify $\proj^1$ with the pencil of lines through $P =
[0,0,1]$, by sending $t$ to the line $z_1 + t \, z_2=0$. Then the four
lines $L, M, PE', PE$ respectively correspond to the $t$-values 
\[ 0, \quad \infty, \quad \pm
\frac{r(r-1)(s^2-s+1)}{s(s-1)(r^2-r+1)}, \] 
leading to a cross-ratio of $-1$. In other words, the line-pair $PE, PE'$ harmonically separates
the line pair $L, M$. This was proved by Rigby in \cite[\S 4]{Rigby} using synthetic methods. 
\subsection{Class field theory} \label{section.class_field.theory} 
In this section, we recall the bare essentials of class field theory
which would suffice to prove Theorems~\ref{theorem.period3}
and~\ref{theorem.period4}. A comprehensive treatment of this subject 
may be found in~\cite{Lang}. The article by Lenstra and
Stevenhagen~\cite{Lenstra} contains an engaging discussion of the Frobenius map
and its relation to the Artin symbol. The computations below were done
in Sage and confirmed in {\sc Maple}. 

Consider the polynomial 
\[ u(x) = x^3-84x^2+896x-1984, \] 
from Section~\ref{section.period3}. 
Let $K$ denote the splitting field of $u(x)$, with discriminant 
$\Delta = 2^{12} \, 13^2 \, 31^2$ and $G = \Gal(K/\QQ) \simeq C_3$. For
every prime number $p$ not dividing $\Delta$, we have an 
Artin symbol $\varphi_p = \artin{K}{\QQ}{p} \in G$. If $u(x)$ factors in $\field_p$
as 
\[ u(x) = u_1(x) \, u_2(x) \dots u_r(x), \] 
then each $u_i$ has the same degree, and moreover it equals
the order of $\varphi_p$ in $G$. Hence $u(x)$ has a root in
$\field_p$, exactly when $\varphi_p$ is the identity element. 

The Kronecker-Weber theorem implies that we 
can embed $K$ in a cyclotomic field $\QQ(\zeta_m)$. Finding such an
$m$ \emph{systematically} is not easy (this would involve
calculating the conductor of the field extension), but we can take
advantage of the fact that the prime factors of $m$ also divide
$\Delta$. If we experiment with small integers made of prime divisors
$2, 13$ and $31$, then we find that $u(x)$ splits
completely\footnote{The article by Roblot~\cite{Roblot} contains a discussion of
  techniques for factoring a polynomial over a number field.} in $L = \QQ(\zeta_{13})$. Thus we have
field extensions $\QQ \subseteq K \subseteq L$, and there is a surjection of
Galois groups 
\[ \Gal(L/\QQ) \stackrel{h}{\longrightarrow} \Gal(K/\QQ), \] 
which sends the Artin symbol $\psi_p = \artin{L}{\QQ}{p}$ to
$\varphi_p$. Hence $\varphi_p=e \iff \psi_p \in \ker h$. 
Now $\Gal(L/\QQ)$ is canonically isomorphic to the
multiplicative group of units 
\[ \ZZ_{13}^\times = \{a \in \ZZ_{13}: (a,13) =1 \} \simeq C_{12}, \] 
and $\psi_p$ acts on $L$ by $z \longrightarrow z^p$. Thus 
\[ \psi_p \in \ker h \iff \psi_p^4=e \iff p^4 \equiv 1
\mod{13}. \] 
The last condition is equivalent to $p \equiv 1, 5, 8, 12
\mod{13}$. 

\subsection{} 
To see this differently, we can express a root of $u(x)$ in terms of
$\zeta = \zeta_{13}$. The actual splitting shows that if we write 
\[ \gamma = 1-\zeta^2 - \zeta^{3} - \zeta^{10} - \zeta^{11} \quad
\text{and} \quad \delta = \zeta^4 + \zeta^6 + \zeta^7 + \zeta^9, \] 
then $\theta = 20 \, \gamma - 4 \, \delta$ is a root of $u(x)$, and 
$K = \QQ(\theta)$. The characteristic property of the Artin symbol
is that 
\[ \varphi_p(z) \equiv z^p \mod{p}, \quad \text{for any $z \in K$};  \] 
that is to say, it imitates the Frobenius map modulo $p$. 
Now $\varphi_p$ is the identity element if and only if it acts
trivially on $\theta$, i.e., exactly when $\theta^p \equiv \theta \mod{p}$. By the freshman's dream, 
\[ \theta^p \equiv 20 \, (1-\zeta^{2p} - \zeta^{3p} - \zeta^{10p} -
\zeta^{11p}) - 4 \, (\zeta^{4p} + \zeta^{6p} + \zeta^{7p} +
\zeta^{9p}), \] 
 and hence this is tantamount to requiring that multiplication by
$p$ should stabilise the sets $\{2,3,10,11\}$ and $\{4,6,7,9\}$
modulo $13$. It is easy to check that this happens exactly when 
$p \equiv 1, 5, 8, 12$. For instance, if $p \equiv 5$, then 
\[ \{2, 3, 10, 11\} \times 5 = \{10, 15, 50, 55 \} \equiv
\{10,2,11,3\} \mod{13}. \] 

\subsection{} \label{section.cases.3and4} 
The other two cases are similar. 
The cubic polynomial $v(x) = x^3 - 48x^2 +656x - 1856$ splits 
completely in $\QQ(\zeta_7)$, and we are reduced to finding the
kernel of the surjection $\ZZ_7^\times \simeq C_6 \longrightarrow
C_3$. These are precisely the elements whose square is the identity, and hence 
\[ \text{$v(x)$ splits in $\field_p$} \iff p^2 
\equiv 1 \mod{7} \iff p \equiv 1, 6 \mod{7}. \] 

The quartic polynomial $w(x) = x^4-144x^3+3152x^2-16896x+25856$ splits completely
in $\QQ(\zeta_{17})$, and we are reduced to finding the kernel of 
$\ZZ_{17}^\times \simeq C_{16} \longrightarrow C_4$. Thus, 
\[ \text{$w(x)$ splits in $\field_p$} \iff p^4 
\equiv 1 \mod{17} \iff p \equiv 1, 4, 13, 16 \mod{17}. \] 

This completes the proofs of Theorems~\ref{theorem.period3} and
\ref{theorem.period4}. \qed 

The criteria in either theorem are inapplicable for finitely many
primes. For instance, in Theorem~\ref{theorem.period3} we must exclude
the primes $2, 13, 31$, as well as those primes $p$ for which $\{f_1, f_2, f_3\}$ is not a Gr{\"o}bner
basis of $I_3$ modulo $p$. It would be possible to cover these
exceptions by a tedious case-by-case computation, but this is unlikely to be of
much interest. 

\subsection{Proof of Theorem~\ref{theorem.2conics}} 
\label{section.conic_computation} 
The theorem follows from some straightforward but elaborate computations,
which were done in {\sc Maple}. Once again, we will write down only
the essential steps. 

Following the definition of the point $P_q(t)$, its coordinates are
found to be 
\[ 
 [ 28 \, q \, t^2-148 \, q \, t + \dots + 60, \; -16 \, q \, t^2+ \dots
 + 64 \, t, \; 14 \, q \, t^2+ \dots + 4 \, t],  \] 
where we have omitted some intermediate terms for brevity. We are searching for a 
conic $\conic_q$ which remains invariant under the map $\tau$. 
Now find the coordinates of $P' = \tau(P)$ and calculate the expression 
\[ f(q,t) = \langle P'R_1, P'R_2, P'R_3,P'R_4 \rangle. \] 
(It is too bulky to be written down in full.) If 
$\tau(\conic_q) = \conic_q$, then we must have $f(q,t) =q$ regardless
of the value of $t$. Hence we solve the equation $\frac{\partial f(q,t)}{\partial t} \equiv 0$. The
derivative factors nicely, and we get a unique solution $q=2/7$. After
back-substitution, we find that $f(2/7,t) =2/7$ as required. Hence
$\tau(\conic_{2/7}) = \conic_{2/7}$. 

Now find the point $Q = \pi(P_{2/7}(t))$, and calculate the cross-ratio $\langle QR_1, QR_2, QR_3,
QR_4 \rangle$. After some simplification, it turns out to be identically equal to 
$-1/4$. This shows that $\pi(\conic_{2/7}) = \conic_{-1/4}$. 
Finally, the formulae for $\imath(t)$ and $\beta(t)$ are determined by the identities 
\[ L_{\imath(t)} = P'R_1, \quad L_{\beta(t)} = QR_1. \] 
This completes the proof. \qed 

As mentioned earlier, the entire calculation is speculative in the
sense that a solution is by no means guaranteed \emph{a priori}. But
fortunately we do get a unique solution with very appealing
properties. 

\subsection{Proof of Theorem~\ref{theorem.leisenring_steiner}}
\label{section.leisen.proof} 
Let $F$ (resp.~$F'$) be the point of intersection of the line triple coming from
even (resp.~odd) permutations. Their coordinates are 
\[ \begin{aligned} 
F & = [r (r-1)(s^2-s+1), s(s-1)(r^2-r+1), (r-1)(s-1)(r+s)(rs-1)], \\ 
F' & = [r (r-1) (s^2-s+1), -s(s-1)(r^2-r+1),
-(r-1)(s-1)(rs+1)(r-s)]. 
\end{aligned} \] 
These points are different from the $E$ and $E'$ of Section~\ref{section.signature.proof}, but essentially the same
argument shows that $F$ and $F'$ are both defined
and distinct, if and only if either $r$ or $s$ is not a root of
$Z^2-Z+1=0$. Let $\Psi_* = FF'$, and now find the cross-ratios
\[ d = \langle \Psi_e, \Psi_{(1 \, 2 \, 3)}, \Psi_{(1 \, 3 
  \, 2)}, \Psi_* \rangle, \qquad 
d' = \langle \Psi_{(1 2)}, \Psi_{(1 \, 3)}, \Psi_{(2 \,
  3)}, \Psi_* \rangle. \] 
A direct calculation shows that $c=d$ and $c'=d'$, hence 
$\leisen(\Pi)$ is equivalent to $\steiner(\Pi)$. \qed 

\begin{figure}[H]
\includegraphics[width=10cm]{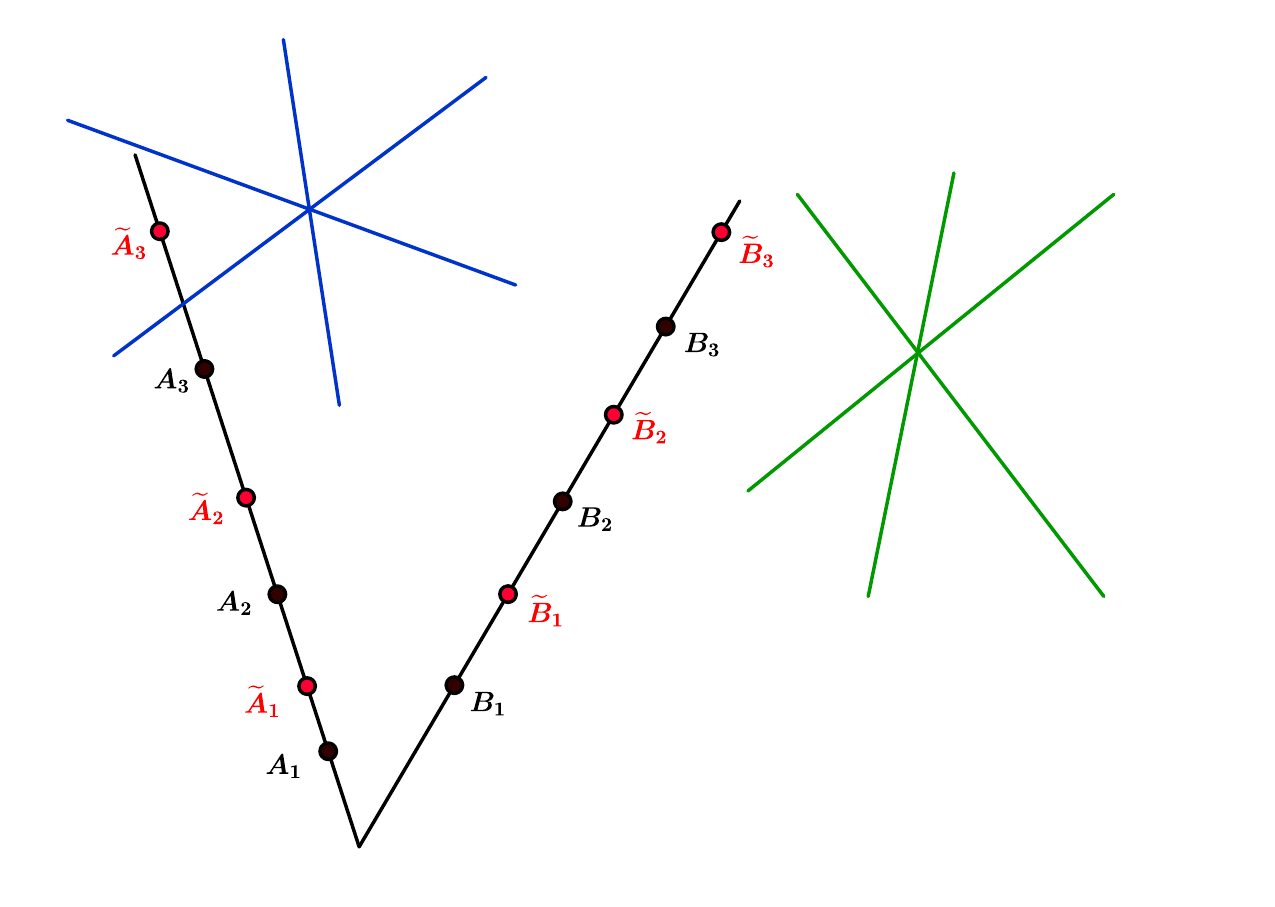} 
\caption{The drawing is not intended to be geometrically accurate. The
actual Pappus lines are not the ones shown.} 
\label{diag:steiner_leisen} 
\end{figure}

Thus we must have an automorphism $\varphi_N$ of $\proj^2$ such that 
$\varphi_N(\Lambda_\sigma) = \Psi_\sigma$ for all $\sigma$ in
$\Sym_3 \cup \{*\}$. Since we know the coordinates of all the fourteen lines, 
we can find its matrix by solving a system of linear
equations. The answer comes
out to be $N = \left[ \begin{array}{ccc} 2 & 0 & s+1 \\ 0 & 2 & r+1 \\
0 & 0 & -1 \end{array} \right]$. Notice that $\varphi_N$
preserves the lines $z_1=0, z_2=0$. Hence, it slides the
points $A_1, \dots, B_3$ along the same lines into new positions
$\widetilde{A}_1, \dots, \widetilde{B}_3$, like beads woven on a
wire, in such a way that 
\[ \text{Pappus lines of $A_1, \dots, B_3$} = 
\text{Leisenring lines of $\widetilde{A}_1, \dots,
  \widetilde{B}_3$}. \] 
This is shown schematically in Diagram~\ref{diag:steiner_leisen}. 

\subsection{} \label{section.involution.pappus} 
Recall the involutive operation $r \rightarrow r^\circ$ from
Section~\ref{section.circ}. Although $\Pi(r,s)$ and $\Pi(r^\circ, s^\circ)$ are in general 
inequivalent Pappus structures, we have 
\[ \sig(\steiner(\Pi(r,s)) = \sig(\steiner(\Pi(r^\circ,s^\circ)). \]
In other words, the two sequences of points $C(r,s)$ and
$C(r^\circ,s^\circ)$ lead  to projectively equivalent sextuples
of Pappus lines. Exactly as in the previous section, we can find the
matrix $T$ whose associated automorphism $\varphi_T$ takes the first
sextuple to the second. The answer comes out to be 
\[ T = \left[ \begin{array}{ccc} 
\frac{(2s-1)(s-2)}{s(s-1)} & 0 & \frac{s-2}{s-1} \\ 
0 & \frac{(2r-1)(r-2)}{r(r-1)} & \frac{r-2}{r-1} \\ 
0 & 0 & -3 
\end{array} \right]. \] 
Once again, it is easy to see that $\varphi_T$ preserves the lines $z_1=0, z_2=0$. 
Said differently, we have two distinct sequences of points, namely $\varphi_T(C(r,s))$ and $C(r^\circ,
s^\circ)$, supported on the same pair of lines, which lead to the 
same sextuple of Pappus lines. It would be interesting to have a synthetic proof of this result. 

\newpage 
\section{Appendix: The Pappus-Steiner map over real numbers} 
\centerline{by \sc{Attila D{\'e}nes}}
This section contains a short computer-aided analysis of the Pappus-Steiner map (\ref{formula.pappus_steiner})
over the field $k = \RR$. In order to study its dynamical behaviour, we used the computer
program developed by D{\'e}nes and Makay~\cite{denes_makay}, which
helps in the visualisation of attractors and basins of dynamical systems. 

\begin{figure}[H]
\includegraphics[height=8cm]{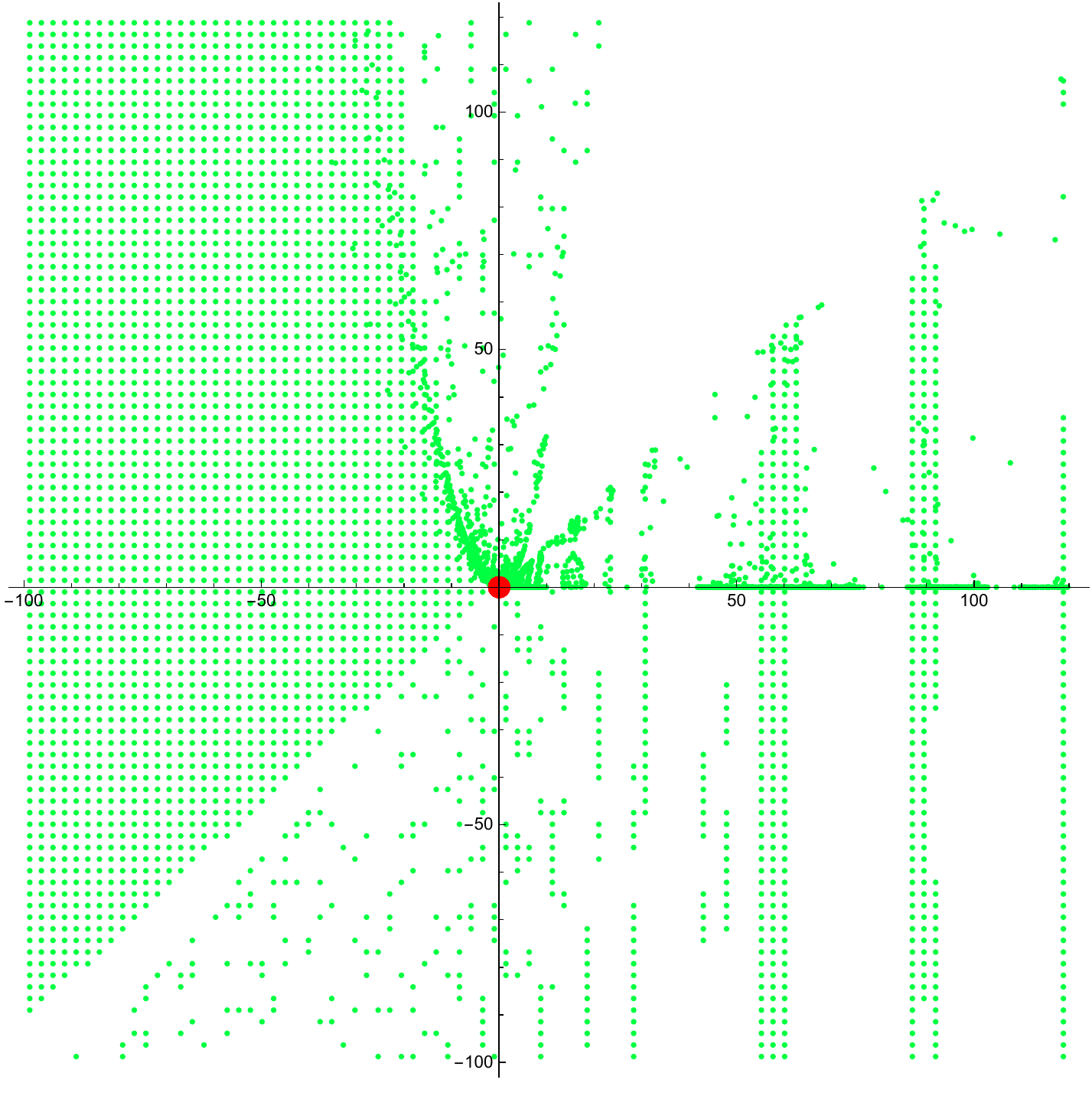}
\caption{The origin as attractor (shown in red)  and its basin (shown
  in green)}
\label{fig:centralpoint}
\end{figure}

According to the algorithm given in~\cite{denes_makay}, if an 
orbit leaves some multiple of the area under examination, then we
assume that this orbit diverges. However, such a judgment must always
be tentative, since such an orbit may return several steps later. If
one calculates some orbits of the Pappus--Steiner map, then one can see
examples of orbits getting very far from the origin and still
returning close to it after some iterations. 
The results obtained suggest that $(0,0)$ is an attractor of the
system (see Diagram~\ref{fig:centralpoint}), although the map itself is
undefined at this point. It is generally difficult to handle these orbits
numerically, since the coordinates of the points become very small
near the origin. 

\begin{figure}[H]
\includegraphics[height=7cm]{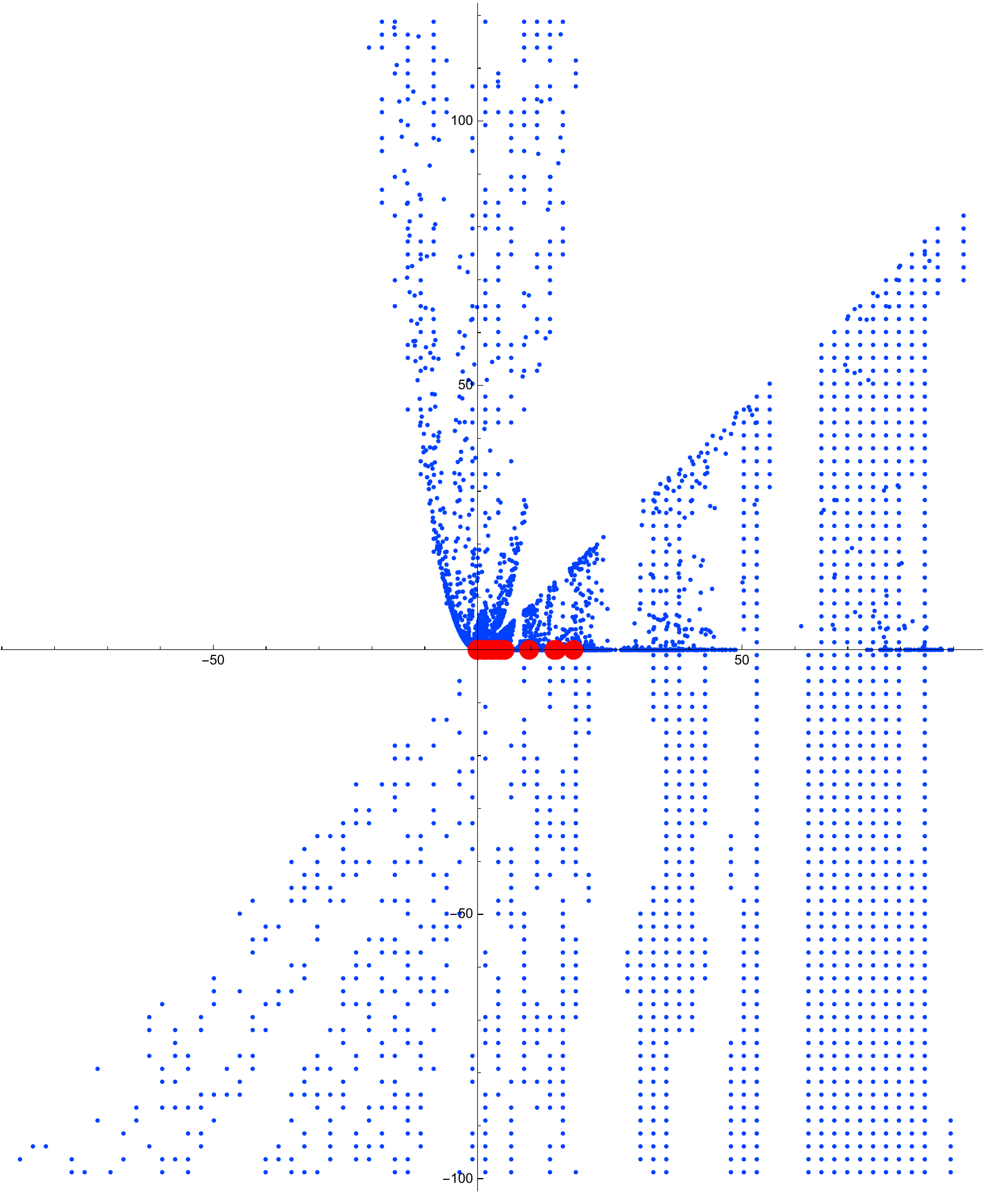}
\caption{Attractor (shown in red) including the origin and points of the $x$ axis and its basin (shown in blue)}\label{fig:attractors}
\end{figure}

There are also orbits which tend to an attractor consisting of the two points $(0,0)$ and $(1,0)$. Although the possibility of numerical errors also holds for these orbits, one may support the conjecture of the existence of such an attractor by considering a point $(1+\varepsilon,\delta)$ close to $(1,0)$, whose first iterate
$$\left(\frac{2\delta(\delta-(1+\varepsilon)+2)}{(1+\varepsilon-\delta)^2},\frac{\delta^2}{(1+\varepsilon-\delta)^2}\right)\approx(2\delta,\delta^2)$$
is close to $(0,0)$, while iterating $(2\delta,\delta^2)$ gives
$$\left(\frac{2\delta^2(\delta^2-2\delta+2)}{(2\delta-\delta^2)^2},\frac{\delta^4}{(2\delta-\delta^2)^2}\right)\approx\left(1-\delta,\frac{\delta^2}{4}\right)$$
which is again close to $(1,0)$. This argument also suggests
asymptotic stability of this attractor, because if $\delta$ is small,
then $\delta^2/4$ is even smaller. One can also obtain similar
attractors formally, e.g.,
$\left\{(0,0),(4,0)\right\}$. Diagram~\ref{fig:attractors} suggests
that there might be other attractors consisting of several points. 

\begin{figure}[H]
\includegraphics[height=7cm]{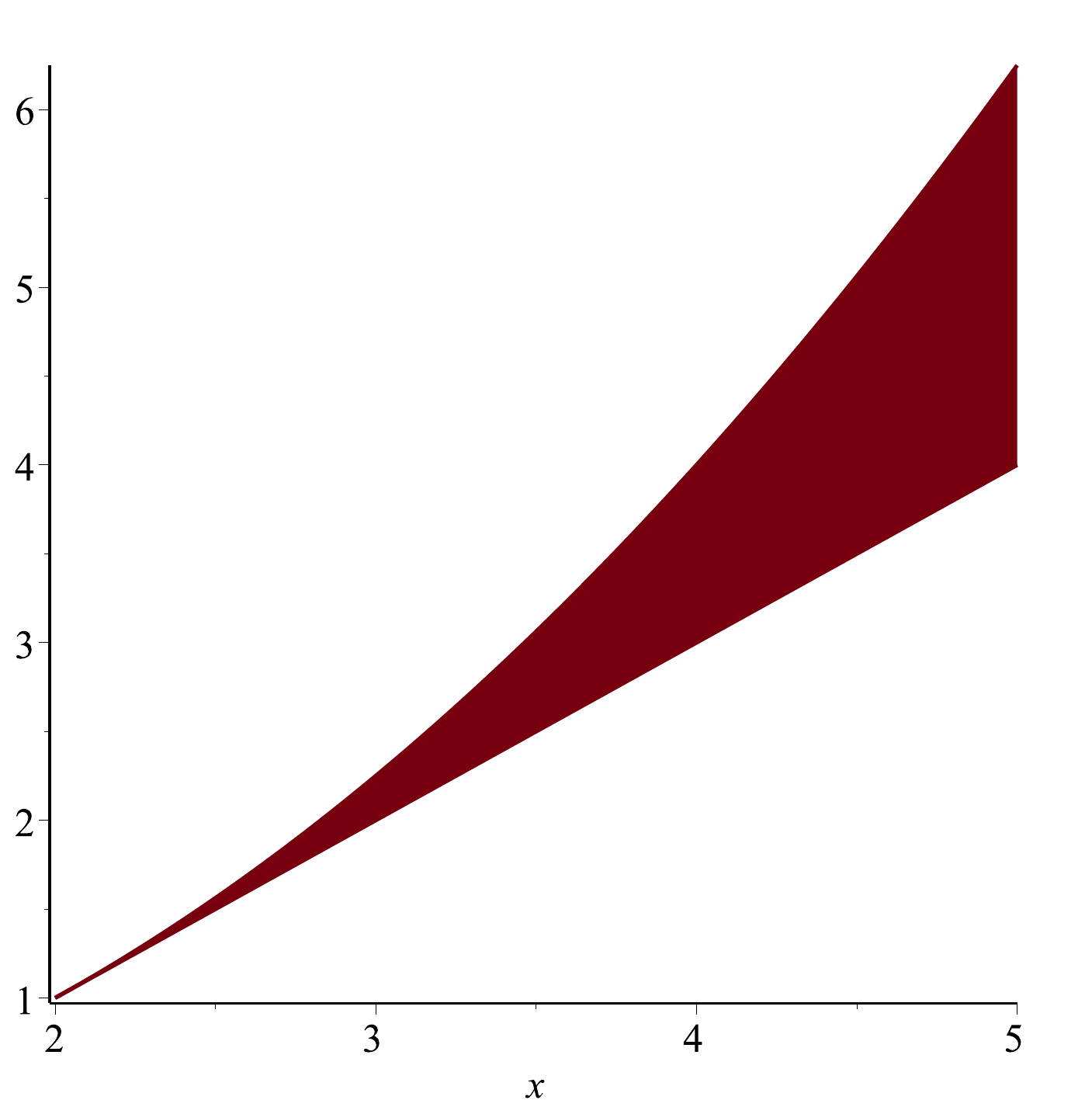} 
\caption{The region between the harmonic line and the balanced
  parabola}
\label{fig:midregion} 
\end{figure}

The most interesting region is the one between the line $y=x-1$ of
harmonic structures, and the parabola $4y = x^2$ of balanced
structures (see Diagram~\ref{fig:midregion}). Numerical experiments suggest
that there are points whose orbits are dense in this region. However 
this statement must again be taken with some degree of caution, since the
expression for $\pi$ becomes numerically unstable near the line $x=y$. 

{\sc Acknowledgement:} The first author is thankful to Siddarth Sankaran for his help 
with the number-theoretic calculations in
Sections~\ref{section.class_field.theory} through~\ref{section.cases.3and4}. 

%%%%%%%%%%%%%%%

\medskip 

\centerline{--} 

\vspace{1cm} 

\parbox{7cm} 
{Jaydeep Chipalkatti \\ 
Department of Mathematics, \\
University of Manitoba, \\ 
Winnipeg, MB R3T 2N2. \\ 
Canada.\\
{\tt jaydeep.chipalkatti@umanitoba.ca}} 
\hfill 
\parbox{7cm} 
{Attila D{\'e}nes \\ 
Bolyai Institute, \\ University of Szeged, \\ 
Aradi v{\'e}rtan{\'u}k tere 1, \\ 
H-6720, Szeged, Hungary. \\ 
{\tt denesa@math.u-szeged.hu}} 

\end{document}